\documentclass{amsart} 
 
\newtheorem{thm}{Theorem}
\newtheorem{lemma}[thm]{Lemma}
\newtheorem{prop}[thm]{Proposition}

\theoremstyle{remark}
\newtheorem{remark}[thm]{Remark}

\newcommand{\End}{\operatorname{End}}
\newcommand{\GK}{\operatorname{GK dim}}
\newcommand{\ann}{\operatorname{ann}}

\begin{document} 
\title{Affine semiprime algebras of GK dimension one are (still) pi}
\author{Christopher J.
Pappacena} \address{Department of Mathematics, Baylor University, Waco, TX
76798} \email{Chris\_\;Pappacena@baylor.edu}
\author{Lance W. Small}
\address{Department of Mathematics, University of California, San Diego, CA
92093}
\email{lwsmall@ucsd.edu}
\author{Jeanne Wald}
\address{Department of Mathematics, Michigan State University, East Lansing, MI 48824}
\email{wald@math.msu.edu}
\subjclass{16P90, 16R20}
\keywords{affine semiprime algebra, GK dimension, polynomial identity}

\thanks{The first author was partially supported by a postdoctoral fellowship from the Mathematical Sciences Research Institute and a Baylor University Summer Sabbatical.}

\begin{abstract} In this note, we give a new proof of the fact
that an affine semiprime algebra $R$ of Gelfand-Kirillov dimension $1$ satisfies a polynomial
identity.  Our proof uses only the growth properties of the algebra and yields an explicit upper bound for the pi degree of $R$.\end{abstract} \maketitle

Let $R$ be an affine algebra over a field $k$.  The Small-Stafford-Warfield Theorem \cite{SSW,SW} asserts that if $R$ has Gelfand-Kirillov (GK) dimension $1$, then $R$ satisfies a polynomial identity.  The proof of this result is in three main steps:  First, the result is proved for $R$ prime in \cite{SW}.  Then, in \cite{SSW} the result is exended to $R$ semiprime by a (somewhat lengthy) proof that $R$ has finitely many minimal primes.  Finally, the result is proven in the general case (again in \cite{SSW}) by showing that the nilradical of $R$ must be nilpotent. 

The purpose of this note is to give a short, self contained proof of the theorem for semiprime algebras.  This new proof is advantageous for several reasons.  First, the proof emphasizes the growth of the algebra and hopefully gives more insight into ``why" the theorem is true.  Second, the argument is less technical; for example, we need only change fields once, instead of three times.  Finally, the proof gives an explicit upper bound for the pi degree of $R$, which is an improvement over the upper bound obtained in \cite{P}.

We have not defined the GK dimension and in fact will not need its definition in the sequel.  Instead we will use the following notion which at first glance is weaker than having GK dimension $1$.  Since $R$ is affine, we can fix a finite set $S$ which generates $R$ as a ring.  Then $R$ is said to have \emph{linear growth} if $R$ is not finite-dimensional over $k$ and $\dim_k kS^i-\dim_kkS^{i-1}\leq c$ for some $c>0$ and all $i$, where $kS^i$ is the span of all words in $S$ of length at most $i$. Note that this condition implies that $\dim_k kS^i\leq ci+1$ for all $i$; however, there may be a smaller constant $d$ such that $\dim_kkS^i\leq di+1$.  The property of having linear growth is independent of the choice of generating space; see \cite[Lemma 1.1]{KL}.  G. Bergman proved that an affine $k$-algebra $R$ has GK dimension $1$ if and only if $R$ has linear growth \cite{Be}, \cite[Theorem 2.5]{KL}. 

Note that the constant $c$ above depends on the choice of generating set $S$.  We shall call the smallest $c$ such that $\dim_k kS^i-\dim_k kS^{i-1}\leq c$ for all $i$ the \emph{Bergman bound} for $S$.  In order to work with an invariant of the ring $R$, we define $c(R)$ by 
\[c(R)=\min\{c:\text{$c$ is the Bergman bound for some generating set $S$}\}.\]

The following result is contained in the original proof of the theorem in the prime case \cite[p. 387]{SW}.  We shall use it below, but feel it is pointless to give a proof because we would essentially be copying the proof given in \cite{SW}.  We remark that the proof of this result uses only the fact that $R$ has linear growth.

\begin{lemma}  If $R$ is a prime affine $k$-algebra of GK dimension $1$, then $R$ satisfies the ascending chain condition on right ideals of the form $\ann r$.\label{SW lemma}
\end{lemma}

\begin{lemma}  Let $k$ be a field, and let $S$ be a finite generating set for
$M_n(k)$ such that
\begin{equation}\dim_k kS^i-\dim_k kS^{i-1}\leq
\ell\label{bound}\end{equation}
for some fixed positive integer $\ell$ and  for
all $i$. Then, any word in $S$ of length $\ell(n+1)$ can be written as a linear
combination of words of strictly smaller length. \label{comb lemma}\end{lemma}

\begin{proof} Enumerate the elements of $S$ as $\{s_1,\dots,s_t\}$, and order $S$ by $s_1<s_2<\dots<s_t$.  We shall denote by $\Sigma(q)$ the set of all sequences $(a_1,\dots, a_q)$, with $a_i\in S$.  We give $\Sigma(q)$ the lexicographic ordering induced by the ordering of $S$. 

Let us set $m=\ell(n+1)$ to ease notation, and let $W$ be the set of all words in $S$ of length $m$ that are not in $kS^{m-1}$.  Assume that $W$ is nonempty, and let $\tilde w=(a_1,\dots,a_m)$ be the smallest element of $\Sigma(m)$ such that $w=a_1\dots a_m\in W$.  Now, $\tilde w$ has exactly $\ell+1$ subsequences of consecutive terms of length $\ell n$; label these as $\tilde v_1,\dots, \tilde v_{\ell+1}$, where $\tilde v_j=(a_j,\dots ,a_{\ell n+j-1})$.  Denote the corresponding subwords of $w$  by $v_1,\dots,v_{\ell+1}$, respectively (so that $v_j=a_j\dots a_{\ell n+j-1}$). 

Suppose first that each of the $v_j$ are distinct elements of $M_n(k)$.  Then, by hypothesis, there is a nontrivial dependence relation $\sum_j \alpha_j v_j+u=0$, with each $\alpha_j\in k$ and $u\in kS^{\ell n-1}$.  Let $\tilde v_p$ be the largest sequence in $\Sigma(\ell n)$ such that $\alpha_p\neq 0$. Then we can solve for $v_p$ as a linear combination of the other terms: $v_p=\sum_{j\neq p}\alpha_p^{-1}\alpha_j v_j +\alpha_p^{-1}u$.  Making this substitution into the word $w$, we have expressed $w$ as a linear combination of words of length $\leq m$.  The key observation is that, by construction, each word of length $m$ in this expression for $w$ has a corresponding sequence in $\Sigma(m)$ which is smaller than $\tilde w$.  Thus, by hypothesis, these terms are not in $W$, and so neither is $w$.

We may therefore suppose that there are two indices $p<q$ such that $v_p=v_q$. We first show that this implies the stronger statement that $\tilde v_p=\tilde v_q$. To see this, we can write 
\begin{equation}w=a_1\dots a_{p-1}(v_p)a_{p+\ell n}\dots a_m=a_1\dots a_{p-1}(v_q)a_{p+\ell n}\dots a_m.\end{equation}
By hypothesis $\tilde w$ is the smallest element of $\Sigma(m)$ representing $w$; consequently, we see that 
\begin{equation}
(a_1,\dots,a_{p-1},a_q,\dots,a_{q+\ell n-1},a_{p+\ell n},\dots, a_m)\geq \tilde w.\end{equation}
This shows that $a_q\geq a_p$.  A similar argument, using 
\begin{equation}
w=a_1\dots a_{q-1}(v_q)a_{q+\ell n}\dots a_m=a_1\dots a_{q-1}(v_p)a_{q+\ell n}\dots w_a\end{equation}
shows that $a_p\geq a_q$, and an induction argument then shows that the two sequences $\tilde v_p$ and $\tilde v_q$ are identical.

Hence, the initial sequences of length $q-p$ of $\tilde v_p$ and $\tilde v_q$ are identical.  Set $\tilde u=(a_p,\dots,a_{q-1})=(a_q,\dots,a_{2q-p-1})$, so that 
\begin{equation}\begin{split}
\tilde v_p&=(a_p,\dots,a_{q-1},a_q,\dots,a_{2q-p-1},a_{2q-p},\dots, a_{p+\ell n-1})\\ &=(\tilde u,\tilde u,a_{2q-p},\dots,a_{p+\ell n-1}).\end{split}\end{equation}
(Note that $q\leq\ell+1\leq \ell n+1$, so that $a_q$ does appear $\tilde v_p$ as claimed.)  Now, an easy induction shows that 
\begin{equation}\tilde v_p=(\tilde u,\dots, \tilde u,a_{nq-(n-1)p},\dots, a_{p+\ell n-1}),\end{equation}
 where $\tilde u$ is repeated $n$ times.  In the extreme case that $q=\ell+1$ and $p=1$, then in fact $\tilde v_p=(\tilde u,\dots,\tilde u)$ (repeated $n$ times) and there are no additional entries in $\tilde v_p$.

If we set $u=a_p\dots a_{q-1}$, then $w=a_1\dots a_{p-1}u^na_{nq-(n-1)p}\dots a_m$. (Again, in the extreme case $q=\ell+1$, $p=1$, we have $w=a_1\dots a_{p-1}u^n$.) Since $w$ contains a subword that is an $n$-th power, we can use the Cayley-Hamilton Theorem to write $w$ as a linear combination of words of length less than $m$. This shows that $W$ is empty, which is precisely the assertion of the lemma.
\end{proof}

Let $\bar k$ denote the algebraic closure of $k$.  We shall call a representation $\rho:R\rightarrow M_n(k)$ \emph{irreducible} if the image of $\rho$ under the inclusion $M_n(k)\subseteq M_n(\bar k)$ generates $M_n(\bar k)$ as an algebra; that is, $\bar k\rho(R)=M_n(\bar k)$.  Our next result bounds the dimensions of the irreducible representations of $R$ by a function of $c(R)$ when $R$ has linear growth.

\begin{prop}  Let $R$ be an affine algebra of GK dimension $1$.  If $\rho:R\rightarrow M_n(k)$ is an irreducible
representation, then $n$ is bounded by a function of $c(R)$.\label{irrep}
\end{prop}

\begin{proof}  Fix a generating set
$T$ for $R$ with Bergman bound $c(R)$, and let
$S=\rho(T)$.  Then, $S$ generates $M_n(\bar k)$ by definition,
and since $\dim_k kT^i-\dim_k kT^{i-1}\leq c(R)$ for all $i$, it follows that
$\dim_{\bar k}\bar k S^i-\dim_{\bar k}\bar k S^{i-1}\leq c(R)$ as well.  By Lemma \ref{comb
lemma}, it follows that every word in $S$ of length $c(R)(n+1)$ is reducible, so that the
words of length $c(R)(n+1)-1$ span $M_n(\bar k)$.  Hence $n^2=\dim_{\bar k} \bar k S^{c(R)(n+1)-1}\leq c(R)^2(n+1)-c(R)+1$.  Solving for $n$ as a function of $c(R)$, we obtain the formula 
\begin{equation}
n\leq {\textstyle\frac{1}{2}}\bigl(c(R)^2+\sqrt{c(R)^4+4c(R)^2-4c(R)+4}\bigr).
\end{equation}
So, the dimensions of the irreducible representations of $R$ are bounded by a constant which depends only on the growth of $R$, as claimed.
\end{proof}

The next step in the proof is to establish that certain primitive $k$-algebras of GK dimension at most $1$ are in fact finite-dimensional over $k$.

\begin{prop}  Let $k$ be an uncountable field and let $R$ be an affine, primitive $k$-algebra with $\GK R\leq 1$.  Then $R$ is finite-dimensional over $k$.  \label{prim}
\end{prop}

\begin{proof}
Fix a generating set $S$ for $R$ with Bergman bound $c$, and let $M$ be a faithful simple $R$-module.  We shall show that $M$ is finite-dimensional over $k$.  It will then follow from standard structure theory that $R$ is also finite-dimensional over $k$.  

So assume to the contrary that $M$ is infinite-dimensional over $k$, and let $x_1,\dots, x_t$ be $D$-linearly independent elements of $M$, where $D=\End_R(M)$.  Note that since $k\subseteq D$, $x_1,\dots,x_t$ are $k$-linearly independent as well.  By the Density Theorem, there exist elements $r_1,\dots, r_t$ such that $x_jr_i=0$ for $j<i$, and $x_ir_i=x_1$.  
Let $A_i=\ann\{x_1,\dots,x_i\}$.  Then, there is an $R$-module surjection
$A_i+r_iR/A_i\rightarrow M$ given by $A_i+r_ir\mapsto x_ir_ir=x_1r$.  (Note that $x_1$ generates $M$ as $M$ is simple.)

Let $N$ be a positive integer such that each of $r_1,\dots, r_t$ are in $kS^N$. The filtration $A_t\leq A_{t-1}\leq\dots\leq A_2\leq A_1$ gives rise to a vector space decomposition
\[A_1\cong A_{1}/A_2\oplus A_{2}/A_{3}\oplus\dots\oplus A_{t-1}/A_t\oplus A_t.\]

This gives rise to the vector space decomposition

\begin{equation}
A_1\cap kS^{n+N}\cong\bigoplus_{i=1}^{t-1} \frac{A_i\cap kS^{n+N}}{A_{i+1}\cap kS^{n+N}}\oplus  A_t\cap kS^{n+N}\label{vec}\end{equation}
for all $n>0$.  Now, since $r_i\in kS^N$, we have that 
\[A_i\cap kS^{n+N}+r_ikS^n\subseteq A_{i-1}\cap kS^{n+N}\] for all $i$. There is a vector space surjection 
\begin{equation}\frac{A_i\cap kS^{n+N}+r_ikS^n}{A_i\cap kS^{n+N}}\rightarrow x_1kS^n\end{equation}
induced by the $R$-module surjection $A_i+r_iR/A_i\rightarrow M$ defined above.  If we set $g(n)=\dim_k x_1kS^n$, then we see that the total dimension of the right hand side of equation \eqref{vec} is at least $(t-1)g(n)$ (one for each of the first $t-1$ terms).  On the other hand, since $A_1\cap kS^{n+N}\subseteq kS^{n+N}$, the total dimension of the left hand side is at most $c(n+N)+1$, which gives the inequality \begin{equation}
c(n+N)+1\geq (t-1)g(n),\label{ineq}\end{equation}
 which holds for all $n\geq 0$.  

Since $M$ is infinite-dimensional over $k$, we must have $g(n)\geq n$: the jump in dimension from $x_1kS^n$ to $x_1kS^{n+1}$ must be at least $1$ unless it stablizes, but it stabilizes if and only if $M$ is finite-dimensional over $k$.  Hence, we deduce from \eqref{ineq} that $c(n+N)+1\geq (t-1)n$, or $(c(n+N)+1)/n\geq t-1$.  Since this inequality holds for arbitrary $n$, we can let $n\rightarrow \infty$ to obtain $c\geq t-1$. We have therefore shown that any set of $D$-linearly independent elements of $M$ has cardinality at most $c+1$.  This shows that $M$ is finite-dimensional over $D$, and that $R\cong M_t(D)$ for some $t\leq c+1$.  Now $D$ is a countable-dimensional $k$-algebra, and $k$ is uncountable, so that $D$ is algebraic over $k$. Thus $\bar k$ is a splitting field for $D$, and so some finite extension $K/k$ splits $D$.  But then $[D:K]$ is finite, so that $[D:k]$ is as well.  Hence $R$ is finite-dimensional over $k$.
\end{proof}

We now have all of the ingredients in place to prove the theorem.  Since we have not yet stated it formally, we do so now.

\begin{thm}  Let $R$ be a semiprime affine algebra of GK dimension $1$.  Then $R$ satisfies a polynomial identity.
\end{thm}

\begin{proof}  Our first step is to reduce to the case that $k$ is uncountable.  Note that if $K/k$ is a field extension, then the canonical map $R\rightarrow R_K:= R\otimes_k K$ is injective.  If we let $K=k(x_i:i\in I)$ where $I$ is an uncountable index set and the $x_i$ are indeterminates over $k$, then $R_K$ is still semiprime, affine over $K$, and of GK dimension $1$.  If we can show that $R_K$ is pi, it will then follow that $R$ is pi.  Thus we replace $R$ with $R_K$ and $k$ with $K$ and (changing notation) we may assume that the base field $k$ is uncountable. 

Since $R$ is semiprime, we can write $R$ as a subdirect product $R\leq \prod R/I$, where each $R/I$ is prime (necessarily of GK dimension $\leq 1$).  By Lemma \ref{SW lemma}, each $R/I$ has the ascending chain condition on right annihilators.  Thus, by \cite[Proposition 2.6.24]{R}, $R/I$ does not have any nil ideals.  Since $k$ is uncountable, \cite[Theorem 2.5.22]{R} shows that $J(R/I)$ is nil, hence $J(R/I)=0$.  We have shown that each $R/I$ is semiprimitive, and since $R$ is a subdirect product of the $R/I$, $R$ is semiprimitive as well.

Thus $R$ embeds in $\prod_{P\in\mathcal{P}} R/P$, where $\mathcal{P}$ is the set of primitive ideals of $R$.  Since each $R/P$ is primitive of GK dimension $\leq 1$, it is finite-dimensional over $k$ by Proposition \ref{prim}.  Moreover, identifying $R/P$ as a subalgebra of $M_n(\bar k)$ for some $n$, we see that the surjection $R\rightarrow R/P$ gives an irreducible representation of $R$.  Thus, $n$ is bounded by the growth of $R$.  In particular, each of the terms in the product $\prod R/P$ has pi degree at most
\begin{equation} N:=\textstyle\frac{1}{2}\bigl(c(R)^2+\sqrt{c(R)^4+4c(R)^2-4c(R)+4}\bigr).\end{equation}
It follows that $R$ has pi degree at most $N$, as well. 
\end{proof}

\begin{remark} (a) The bound of $N$ for the pi degree of $R$ is roughly half the bound obtained in \cite{P}. 

(b) An easy to remember upper bound for the pi degree of $R$ is $c(R)^2+1$. 
\end{remark}

\bibliographystyle{amsalpha}

\end{document}